\documentclass{gtpart}
\usepackage{amsmath,amssymb,amsthm,stmaryrd}
\usepackage[all]{xy}
\title{An overview of abelian varieties in homotopy theory}
\author{Tyler Lawson}
\givenname{Tyler}
\surname{Lawson}
\address{Department of Mathematics\\University of Minnesota\\Minneapolis, MN 55455\\USA}
\email{tlawson@math.umn.edu}


\newtheorem{thm}{Theorem}

\theoremstyle{definition}

\theoremstyle{remark}

\newtheorem{exam}[thm]{Example}

\newtheorem{slogan}[thm]{Slogan}
\newtheorem{ques}[thm]{Question}

\def\co{\colon\thinspace}
\newcommand{\mb}[1]{\mathbb{#1}}
\newcommand{\mf}[1]{\mathfrak{#1}}

\newcommand{\Ext}{\ensuremath{{\rm Ext}}}

\newcommand{\overto}{\mathop\rightarrow}

\newcommand{\onto}{\mathop\twoheadrightarrow}
\newcommand{\longoverto}{\mathop{\longrightarrow}}

\newcommand{\Uni}{{\rm U}}
\newcommand{\BU}{{\rm BU}}

\newcommand{\TAF}{{\rm TAF}}

\newcommand{\End}{{\rm End}}

\newcommand{\As}{{\cal A}s}
\newcommand{\LT}{{\rm LT}}
\newcommand{\Spec}{{\rm Spec}}

\newcommand{\tmf}{{\rm tmf}}
\newcommand{\eo}{{\rm eo}}
\newcommand{\TMF}{{\rm TMF}}
\newcommand{\Ell}{{\cal E}ll}

\newcommand{\Sh}{{\rm Sh}}

\newcommand{\eilm}[1]{\ensuremath{{\mb H} #1}}
\newcommand{\smsh}[1]{\ensuremath{\mathop{\wedge}_{#1}}}

\newcommand{\form}[2]{\ensuremath{\left\langle#1,#2\right\rangle}}

\newcommand{\xym}[1]{
\vskip 0.7pc
\centerline{\xymatrix{#1}}
\vskip 0.7pc
}

\begin{document}

\begin{abstract}
  We give an overview of the theory of formal group laws in homotopy
  theory, leading to the connection with higher-dimensional abelian
  varieties and automorphic forms.
\end{abstract}

\maketitle

\section{Introduction}

The goal of this paper is to provide an overview of joint work with
Behrens on topological automorphic forms \cite{taf}.  The ultimate
hope is to introduce a somewhat broad audience of topologists to this
subject matter connecting modern homotopy theory, algebraic geometry,
and number theory.

Through an investigation of properties of Chern classes, Quillen
discovered a connection between stable homotopy theory and
$1$-dimensional formal group laws \cite{quillen}.  After almost 40
years, the impacts of this connection are still being felt.  The
stratification of formal group laws in finite characteristic gives
rise to the {\em chromatic filtration} in stable homotopy theory
\cite{ravenelgreen}, and has definite calculational consequences.  The
nilpotence and periodicity phenomena in stable homotopy groups of
spheres arise from a deep investigation of this connection \cite{dhs}.

Formal group laws have at least one other major manifestation: the
study of abelian varieties.  The examination of this connection led to
elliptic cohomology theories and topological modular forms, or $\tmf$
\cite{hopkinsicm}.  One of the main results in this theory is the
construction of a spectrum $\tmf$, a structured ring object in the
stable homotopy category.  The homotopy groups of $\tmf$ are, up to
finite kernel and cokernel, the ring of integral modular forms
\cite{deligne} via a natural comparison map.  The spectrum $\tmf$ is
often viewed as a ``universal'' elliptic cohomology theory
corresponding to the moduli of elliptic curves.  Unfortunately, the
major involved parties have not yet published a full exposition of
this theory.  The near-future reader is urged to consult
\cite{marktmf}, as well as seek out some of the unpublished literature
and reading lists on topological modular forms if more background
study is desired.

Algebraic topology is explicitly tied to $1$-dimensional formal group
laws, and so the formal group laws of higher-dimensional abelian
varieties (and larger possible ``height'' invariants of those) are
initially not connected to topology.  The goal of \cite{taf} was to create
generalizations of the theory of topological modular forms, through
certain moduli of abelian varieties with extra data specifying
$1$-dimensional {\em summands} of their formal group laws.

The author doubts that it is possible to cover all of this
background to any degree of detail within the confines of a paper of
reasonable size, even restricting to those subjects that are of
interest from a topological point of view.  In addition, there are
existing (and better) sources for this material.  Therefore, our
presentation of this material is informal, and we will try to list
references for those who find some subject of interest to them.  We
assume a basic understanding of stable homotopy theory, and an
inevitable aspect of the theory is that we require more and more of
the language of algebraic geometry as we proceed.

A rough outline of the topics covered follows.

In sections \ref{sec:fgl} and \ref{sec:quillen} we begin with some
background on the connection between the theory of complex bordism and
formal group laws.  We next discuss in section \ref{sec:stack} the
basic theories of Hopf algebroids and stacks, and the relation between
stack cohomology and the Adams-Novikov spectral sequence in section
\ref{sec:anss}.  We then discuss the problem of realizing formal group
law data by spectra, such as is achieved by the Landweber exact
functor theorem and the Goerss-Hopkins-Miller theorem, in section
\ref{sec:realize}.  Examples of multiplicative group laws are
discussed in section \ref{sec:multiplicative}, and the theories of
elliptic cohomology and topological modular forms in sections
\ref{sec:elliptic} and \ref{sec:tmf}.  We then discuss the possibility
of moving forward from these known examples in section
\ref{sec:chromatic}, by discussing some of the geometry of the moduli
of formal groups and height invariants.

The generalization of the Goerss-Hopkins-Miller theorem due to Lurie,
without which the subject of topological automorphic forms would be
pure speculation, is introduced in section \ref{sec:lurie}.  We view
it as our point of entry: given this theorem, what kinds of new
structures in homotopy theory can we produce?

The answer, in the form of various moduli of higher-dimensional
abelian varieties, appears in section \ref{sec:pel}.  Though the
definitions of these moduli are lifted almost directly from the study
of automorphic forms, we attempt in sections \ref{sec:end},
\ref{sec:pol}, and \ref{sec:level} to indicate why this data is {\em
  natural} to require in order produce moduli satisfying the
hypotheses of Lurie's theorem.  In section \ref{sec:questions}, we try
to indicate why some initial choices are made the way they are.

One of the applications in mind has been the construction of finite
resolutions of the $K(n)$-local sphere.  Henn has given finite length
algebraic resolutions allowing computation of the cohomology of the
Morava stabilizer group in terms of the cohomology of finite subgroups
\cite{henn}.  Goerss-Henn-Mahowald-Rezk \cite{ghmr} and Behrens
\cite{marksphere} gave analogous constructions of the $K(2)$-local
sphere at the prime $3$ out of a finite number of spectra of the form
$E_2^{hG}$, where $E_2$ is a Lubin-Tate spectrum and $G$ is a finite
subgroup of the Morava stabilizer group.  The hope is that these
constructions will generalize to other primes and higher height by
considering diagrams of abelian varieties and isogenies.

None of the (correct) material in this paper is new.

\section{Generalized cohomology and formal group laws}
\label{sec:fgl}

Associated to a generalized cohomology theory $E$ with (graded)
commutative multiplication, we can ask whether there is a reasonable
theory of Chern classes for complex vector bundles.

The base case is that of line bundles, which we view as being
represented by homotopy classes of maps $X \to \BU(1) =
\mb{CP}^\infty$ for $X$ a finite CW-complex.  An {\em orientation} of
$E$ is essentially a first Chern class for line bundles.  More
specifically, it is an element $u \in E^2({\mb CP}^\infty)$ whose
restriction to $E^2(\mb{CP}^1) \cong E_0$ is the identity element
$1$ of the ring $E_*$.  For any line bundle $L$ on $X$ represented by
a map $f\co X \to \mb{CP}^\infty$, we have an $E$-cohomology element
$c_1(L) = f^*(u) \in E^2(X)$ which is the desired first Chern class.

Orientations do not necessarily exist; for instance, real K-theory
$KO$ does not have an orientation.  When orientations do exist, we say
that the cohomology theory is {\em complex orientable}.  An
orientation is not necessarily unique; given any orientation $u$, any
power series $v = \sum b_i u^{i+1}$ with $b_i \in E_{2i}, b_0 = 1$
determines another orientation and another Chern class.  Any other
orientation determines and is determined uniquely by such a power
series.

Given an orientation of $E$, we can derive computations of
$E^*(\BU(n))$ for all $n \geq 0$, and conclude that for a vector
bundle $\xi$ on a finite complex $X$ there are higher Chern classes
$c_i(\xi) \in E^{2i}(X)$ satisfying naturality, the Cartan formula,
the splitting principle, and almost all of the desirable properties of
Chern classes in ordinary cohomology.  See \cite{adams}.

The one aspect of this theory that differs from ordinary cohomology
has to do with tensor products.  For line bundles $L_1$ and $L_2$,
there is a tensor product line bundle $L_1 \otimes L_2$ formed by
taking fiberwise tensor products.  On classifying spaces, if $L_i$ are
classified by maps $f_i\co X \to \BU(1)$, the tensor product is
classified by $\mu \circ (f_1 \times f_2)$, where $\mu\co \BU(1)
\times \BU(1) \to \BU(1)$ comes from the multiplication map on
$\Uni(1)$.

There is a universal formula for the tensor product of two line
bundles in $E$-cohomology, given by the formula
\[
c_1(L \otimes L') = \sum a_{i,j} c_1(L)^i c_1(L')^j
\]
for $a_{i,j} \in E_{2i + 2j - 2}$.  This formula is valid for all line
bundles but the coefficients $a_{i,j}$ depend only on the orientation.
We often denote this power series in the alternate forms
\[
\sum a_{i,j} x^i y^j = F(x,y) = x +_F y.
\]
This last piece of notation is justified as follows.  The tensor
product of line bundles is associative, commutative, and unital up to
natural isomorphism, and so by extension the same is true for the
power series $x +_F y$:
\begin{itemize}
\item $x +_F 0 = x$,
\item $x +_F y = y +_F x$, and
\item $(x +_F y) +_F z = x +_F (y +_F z)$.
\end{itemize}
These can be written out in formulas in terms of the coefficients
$a_{i,j}$, but the third is difficult to express in closed form.  A
power series with coefficients in a ring $R$ satisfying the above
identities is called a (commutative, $1$-dimensional) {\em formal
  group law over $R$}, or just a formal group law.

The formal group law associated to $E$ depends on the choice of
orientation.  However, associated to a different orientation $v =
g(u)$, the formal group law $G(x,y) = x +_G y$ satisfies
\[
g(x +_F y) = g(x) +_G g(y).
\]
We say that two formal group laws differing by such a
change-of-coordinates for a power series $g(x) = x + b_1 x^2 + \cdots$
are {\em strictly isomorphic}.  (If we forget which orientation we
have chosen, we have a formal group law without a choice of coordinate
on it, or a {\em formal group}.)

The formal group detects so much intricate information about the
cohomology theory $E$ that it is well beyond the scope of this
document to explore it well \cite{ravenelgreen}.  For certain
cohomology theories $E$ (such as Landweber exact theories discussed in
section \ref{sec:realize}), the formal group determines the cohomology
theory completely.  One can then ask, for some spaces $X$, to
understand the cohomology groups $E^*(X)$ in terms of the formal group
data.  For example, if $X = \BU\langle 6\rangle$, this turns out to be
related to cubical structures \cite{ahs}.

\section{Quillen's theorem}
\label{sec:quillen}

There is a cohomology theory $MU$ associated to complex bordism that
comes equipped with an orientation $u$.  There is also a ``smash
product'' cohomology theory $MU \smsh{} MU$ coming equipped with two
orientations $u$ and $v$, one per factor of $MU$, and hence with two
formal group laws with a strict isomorphism $g$ between them.

The ring $L = MU_*$ forming the ground ring for complex bordism was
calculated by Milnor \cite{milnor}, and similarly for $W = (MU \smsh{}
MU)_*$.  Both are infinite polynomial algebras over $\mb Z$, the
former on generators $x_i$ in degree $2i$, the latter on the $x_i$ and
additional generators $b_i$ (also in degree $2i$).  The following
theorem, however, provides a more intrinsic description of these
rings.

\begin{thm}[Quillen]
  The ring $L$ is a classifying object for formal group laws in the
  category of rings, i.e. associated to a ring $R$ with formal group
  law $F$, there is a unique ring map $\phi\co L \to R$ such that the
  image of the formal group law in $L$ is $F$.

  The ring $W \cong L[b_1,b_2,\ldots]$ is a classifying object for
  pairs of strictly isomorphic formal group laws in the category of
  rings, i.e. associated to a ring $R$ with a strict isomorphism $g$
  between formal group laws $F$ and $G$, there is a unique ring map
  $\phi\co W \to R$ such that the image of the strict isomorphism in
  $W$ is the strict isomorphism in $R$.
\end{thm}

(It is typical to view these rings as geometric objects $\Spec(L)$ and
$\Spec(W)$, which reverses the variance; in schemes, these are
classifying objects for group scheme structures on a formal affine
scheme $\hat{\mb A}^1$.)

The structure of the ring $L$ was originally determined by Lazard, and
it is therefore referred to as the Lazard ring.

There are numerous consequences of Quillen's theorem.  For a general
multiplicative cohomology theory $R$, the theory $MU \smsh{} R$
inherits the orientation $u$, and hence a formal group law.  The
cohomology theory $MU \smsh{} MU \smsh{} R$ has two orientations
arising from the orientations of each factor, and these two differ by
a given strict isomorphism.  For more smash factors, this pattern
repeats.  Philosophically, we have a ring $MU_* R$ with formal group
law, together with a compatible action of the group of strict
isomorphisms.

Morava's survey \cite{moravasurvey} is highly recommended.

\section{Hopf algebroids and stacks}
\label{sec:stack}

The pair $(MU, MU \smsh{} MU)$ and the associated rings $(L,W)$
have various structure maps connecting them.  Geometrically, we
have the following maps of schemes.
\xym{
\Spec(L) \ar[r] & 
\Spec(W) \ar@/_4pt/[l] \ar@/^4pt/[l] \ar@(ul,ur)[]&
\Spec(W) {\times}_{\Spec(L)} \Spec(W) \ar[l]
}
These maps and their relationships are most concisely stated by saying
that the result is a groupoid object in schemes.  We view $\Spec(L)$
as the ``object'' scheme and $\Spec(W)$ as the ``morphism'' scheme,
and the maps between them associate:
\begin{itemize}
\item an identity morphism to each object,
\item source and target objects to each morphism,
\item an inverse to each morphism, and
\item a composition to each pair of morphisms where the source of the
  first is the target of the second.
\end{itemize}
The standard categorical identities (unitality, associativity) become
expressed as identities which the morphisms of schemes must satisfy.

A pair of rings $(A,\Gamma)$ with such structural morphisms is a
representing object for a covariant functor from rings to groupoids;
such an object is generally referred to as a {\em Hopf algebroid}
\cite[Appendix A]{ravenelgreen}.

\begin{exam}
Associated to a map of rings $R \to S$, we have the Hopf algebroid
$(S, S \otimes_R S)$, sometimes called the {\em descent} Hopf
algebroid associated to this map of rings.  This represents the
functor on rings which takes a ring $T$ to category whose objects are
morphisms from $S \to T$ (or $T$-points of $\Spec(S)$), and where two
objects are isomorphic by a unique isomorphism if and only if they
have the same restriction to $R \to T$.

More scheme-theoretically, given a map $Y \to X$ of schemes, we get a
groupoid object $(Y, Y \times_X Y)$ in schemes with the same properties.
\end{exam}

\begin{exam}
If $S$ is a ring with an action of a finite group $G$, then there is a
Hopf algebroid $(S, \prod_G S)$ representing a category of points of
$\Spec(S)$ and morphisms the action of $G$ by precomposition.

Again in terms of schemes, associated to a scheme $Y$ with a (general)
group $G$ acting, we get a groupoid object $(Y, \coprod_G Y)$ in
schemes.  It is a minor but perpetual annoyance that infinite products
of rings do not correspond to infinite coproducts of schemes;
$\Spec(R)$ is always quasi-compact.
\end{exam}

\begin{exam}
If $(A,\Gamma)$ is a Hopf algebroid and $A \to B$ is a map of rings,
then there is an induced Hopf algebroid $(B, B \otimes_A \Gamma
\otimes_A B)$.\footnote{Note that the ``descent'' Hopf algebroid is
a special case.}  The natural map
\[
(A,\Gamma) \to (B, B \otimes_A \Gamma \otimes_A B)
\]
represents a fully faithful functor between groupoids, with the map on
objects being the map from points of $\Spec(B)$ to points of
$\Spec(A)$.  This is an equivalence of categories on $T$-points if and
only if this map of categories is essentially surjective (every
object is isomorphic to an object in the image).

In schemes, if $(X,Y)$ is a groupoid object in schemes and $Z \to X$
is a morphism, there is the associated pullback groupoid $(Z, Z
\times_X Y \times_X Z)$ with a map to $(X,Y)$.
\end{exam}

In principle, for a groupoid object $(X,Y)$ there is an associated
``quotient object,'' the coequalizer of the source and target
morphisms $Y \to X$.  This categorical coequalizer, however, is
generally a very coarse object.  The theories of orbifolds and stacks are
designed to create ``gentle'' quotients of these objects by
remembering {\em how} these points have been identified rather than just
remembering the identification.

To give a more precise definition of stacks, one needs to discuss
Grothendieck topologies.  A Grothendieck topology gives a criterion
for a family of maps $\{U_\alpha \to X\}$ to be a ``cover'' of $X$; for
convenience we will instead regard this as a criterion for a single
map $\coprod U_\alpha \to X$ to be a cover.  The category of stacks in
this Grothendieck topology has the following properties.
\begin{itemize}
\item Stacks, like groupoids, form a 2-category (having morphisms and
  natural transformations between morphisms).
\item The category of stacks is closed under basic constructions such
  as 2-categorical limits and colimits.
\item Associated to a groupoid object $(X,Y)$, there is a functorial
  associated stack $\As(X,Y)$.
\item If $Z \to X$ is a cover in the Grothendieck topology, then the
  map of groupoids $(Z, Z \times_X Y \times_X Z) \to (X,Y)$ induces an
  equivalence on associated stacks.
\end{itemize}
In some sense stacks are characterized by these properties
\cite{hollander}.  In particular, to construct a map from a scheme $V$
to the associated stack $\As(X,Y)$ is the same as to find a cover $U
\to V$ and a map from the descent object $(U, U \times_V U)$ to
$(X,Y)$, modulo a notion of natural equivalence.

Stacks appear frequently when classifying families of objects over a
base.  In particular, in the case of the Hopf algebroid of formal
group laws $(\Spec(L), \Spec(W))$ classifying formal group laws and
strict isomorphisms, the associated stack ${\cal M}_{sFG}$ is
referred to as the {\em moduli stack of formal groups} (and strict
isomorphisms)\footnote{As $L$ and $W$ are graded rings, this moduli
stack inherits some graded aspect as well that can be confusing from
a geometric point of view.  It is common to replace $MU$ with a
2-periodic spectrum $MP$ to remove all gradings from the picture;
the resulting Hopf algebroid arising from $MP$ and $MP \smsh{} MP$
classifies formal group laws and {\em non}-strict isomorphisms, but
has the gradings removed.  The associated stack is usually written
${\cal M}_{FG}$, and has the same cohomology.}.

The theory of stacks deserves much better treatment than this, and the
reader should consult other references
\cite{prsomsfgl,niko,coctalos,vistoli,lmb}.  What this rough outline is
meant to do is perhaps provide some intuition.  Stacks form some
family of categorical objects including quotients by group actions,
having good notions of gluing.  A Hopf algebroid gives a {\em
  presentation}, or a coordinate chart, on a stack.

When algebraic topology studies these topics, it is typically grounded
in the study of Hopf algebroids; the more geometric language of stacks
is adopted more recently and less often.  There are several reasons for
this.

This link to algebraic geometry historically {\em only} occurred
through Hopf algebroids.  The development of structured categories of
spectra has made some of these links more clear, but there is still
some foundational work to be done on coalgebras and comodules in
spectra.

Additionally, the theory and language of stacks are not part of the
typical upbringing of topologists, and have a reputation for being
difficult to learn.  By contrast, Hopf algebroids and comodules admit
much more compact descriptions.

Finally, there is the aspect of computation.  Algebraic topologists
need to compute the cohomology of the stacks that they study, and Hopf
algebroids provide very effective libraries of methods for this.  In
this respect, we behave much like physicists, who become intricately
acquainted with particular methods of computation and coordinate
charts for doing so, rather than regularly taking the ``global''
viewpoint of algebraic geometry.  (The irony of this situation is
inescapable.)

By default, when we speak about stacks in this paper our underlying
Grothendieck topology is the ``fpqc'' (faithfully flat, quasi-compact)
topology.  Most other Grothendieck topologies in common usage are not
geared to handle infinite polynomial algebras such as the Lazard ring.

\section{Cohomology and the Adams-Novikov spectral sequence}
\label{sec:anss}

We fix a Hopf algebroid $(A,\Gamma)$, and assume $\Gamma$ is a flat
$A$-module (equivalently under either the source or target morphism).
We regard the source and target morphisms $A \to \Gamma$ as right and
left module structures respectively.

A {\em comodule} over this Hopf algebroid is a left $A$-module $M$
together with a map of left $A$-modules
\[
\phi\co M \to \Gamma \otimes_A M.
\]
We require that the composite
\[
M \overto^\phi \Gamma \otimes_A M \longoverto^{\epsilon \otimes 1} A \otimes_A M
\]
is the identity, where $\epsilon$ is the augmentation $\Gamma \to A$,
and that the two composites
\[
(c \otimes 1)\phi, (1 \otimes \phi)\phi \co
M \to \Gamma \otimes_A
\Gamma \otimes_A M
\]
are equal, where $c$ is the comultiplication $\Gamma \to \Gamma
\otimes_A \Gamma$.  (This map is typically referred to as a {\em
  coaction} which is counital and coassociative.)

The structure of a comodule is equivalent to having an isomorphism of
$\Gamma$-modules
\[
\Gamma \otimes_A^s M \to \Gamma \otimes_A^t M,
\]
tensor product along the source and target $A$-module structures on
$\Gamma$ respectively, satisfying some associativity typically
appearing in the study of descent data.

The category of $(A,\Gamma)$ comodules forms an abelian category.
This category is the category of {\em quasicoherent sheaves} on the
associated stack ${\cal M} = \As(\Spec(A),\Spec(\Gamma))$.  In
general, one needs to show that homological algebra in this category
can reasonably be carried out; see \cite{franke} for details.

Ignoring the fine details, one can define the {\em coherent
  cohomology} of the stack with coefficients in a comodule $M$ to be
\[
\Ext^*_{\text{q.-c.}/{\cal M}}(A,M) = \Ext^*_{(A,\Gamma)}(A,M).
\]
This is computed by the cobar complex
\[
0 \to M \to \Gamma \otimes_A M \to \Gamma \otimes_A \Gamma \otimes_A M
\to \cdots,
\]
where the boundary maps are alternating sums of unit maps,
comultiplications, and the coaction on $M$.  A better definition is
that these groups are the derived functors of the global section
functor on the stack.  As such, this is genuinely an invariant of the
stack itself, and this underlies many change-of-rings isomorphisms:
for example, for a faithfully flat map $A \to B$ the associated cobar
complex for the comodule $B \otimes_A M$ over $(B, B \otimes_A \Gamma
\otimes_A B)$ computes the same cohomology.  (This is both an
important aspect of the theory of ``faithfully flat descent'' and a
useful computational tactic.)

The importance of coherent cohomology for homotopy theory is the
Adams-Novikov spectral sequence.  For a spectrum $X$, the
$MU$-homology $MU_* X$ inherits the structure of an $(L,W)$-comodule,
and we have the following result.

\begin{thm}
There exists a (bigraded) spectral sequence with $E_2$-term
\[
\Ext^{**}_{(L,W)} (L, MU_* X)
\]
whose abutment is $\pi_* X$.  If $X$ is connective, the spectral
sequence is strongly convergent.
\end{thm}

This spectral sequence arises through a purely formal construction in
the stable homotopy category, and does not rely on any stack-theoretic
constructions.  It is a generalization of the Adams spectral sequence,
which is often stated using cohomology and has $E_2$-term $\Ext$ over
the mod-$p$ Steenrod algebra.

We can recast this in terms of stacks.  Any spectrum $X$ produces a
quasicoherent sheaf on the moduli stack of formal group laws, and
there is a spectral sequence converging from the cohomology of the
stack with coefficients in this sheaf to the homotopy of $X$.  Because
in this way we see ourselves ``recovering $X$ from the quasicoherent
sheaf,'' we find ourselves in the position to state the following.

\begin{slogan}
The stable homotopy category is approximately the category of
quasicoherent sheaves on the moduli stack of formal groups ${\cal
M}_{sFG}$.\footnote{Strictly speaking, one should phrase this in
terms of $MU$-local spectra, which are the only spectra that $MU$
can recover full information about.  The current popular
techniques concentrate on $MU$-local spectra, as they include most
of the examples of current interest and we have very few tactics
available to handle the rest.}
\end{slogan}

This approximation, however, is purely in terms of algebra and it
does not genuinely recover the stable homotopy category.  (The Mahowald
uncertainty principle claims that any algebraic approximation to
stable homotopy theory must be infinitely far from correct.)  However,
the reader is invited to consider the following justification for the
slogan.

An object in the stable homotopy category is generally considered as
being ``approximated'' by its homotopy groups; they provide the basic
information about the spectrum, but they are connected together by a
host of $k$-invariants that form the deeper structure.

The spectrum $MU$ is a highly structured ring object, and the
pair $(MU, MU \smsh{} MU)$ forms a ``Hopf algebroid'' in spectra.  A
general spectrum $X$ gives rise to a comodule $MU \smsh{} X$, and
there is a natural map
\[
X \to F_{(MU,MU \smsh{} MU)}(MU, MU \smsh{} X)
\]
from $X$ to the function spectrum of comodule maps; if we believe in
flat descent in the category of spectra, this map should be a weak
equivalence when $X$ is ``good.''  The Adams-Novikov spectral sequence
would then simply be an algebraic attempt to recover the homotopy of the
right-hand side by a universal coefficient spectral sequence ($\Ext$
on homotopy groups approximates homotopy groups of mapping spaces).

The author is hopeful that the theory of comodules in spectra will
soon be fleshed out rigorously.

We note that, in line with this slogan, Franke has proven that for
$2(p-1) > n^2 + n$, the homotopy category of $E_n$-local spectra at
the prime $p$ is the derived category of an abelian category
\cite{franke}, generalizing a result of Bousfield for $n = 1$
\cite{bousfield}.  As is standard, this excludes the primes where
significant nontrivial behavior is present in the homotopy category.

\section{Realization problems}
\label{sec:realize}

Given our current state of knowledge, it becomes reasonable to ask
questions about our ability to construct spectra.

\begin{enumerate}
\item Can we realize formal group laws by spectra?
\item Can we realize them functorially?
\end{enumerate}

More precisely.

\begin{enumerate}
\item Suppose we have a graded ring $R$ with formal group law $F$.
  When can we construct an oriented ring spectrum $E$ whose homotopy
  is $R$ and whose associated generalized cohomology theory has formal
  group law $F$?
\item Suppose we have a diagram of graded rings $R_\alpha$ and formal
  group laws $F_\alpha$ equipped with strict isomorphisms $\gamma_f\co 
  F_\beta \to f^*F_\alpha$ of formal group laws for any map $f\co
  R_\alpha \to R_\beta$ in the diagram, satisfying $\gamma_g \circ g^*(\gamma_f)
  = \gamma_{gf}$.  When can we realize this as a diagram
  $\{E_\alpha\}$ of ring spectra?
\end{enumerate}

More refined versions of these questions can also be asked; we can ask
for the realizations to come equipped with highly structured
multiplication in some fashion.

Two of the major results in this direction are the Landweber exact
functor theorem and the Goerss-Hopkins-Miller theorem.

We recall \cite[Appendix 2]{ravenelgreen} that for any prime $p$, there is a
sequence of elements $(p,v_1,v_2,\ldots)$ of $L$ such that, if $F$ is the
universal formal group law over the Lazard ring $L$,
\[
[p](x) = x +_F \cdots +_F x \equiv v_n x^{p^n} \mod (p,v_1,\cdots,v_{n-1}).
\]
The elements $v_n$ are well-defined modulo lower elements, but there
are multiple choices of lifts of them to $L$ (such as the Hazewinkel
or Araki elements) that each have their advocates.  (By convention,
$v_0 = p$.)

Associated to a formal group law over a field $k$ classified by a map
$\phi\co L \to k$, there are {\em height} invariants
\[
ht_p(F) = \inf\{n\ |\ \phi(v_n) \neq 0\}
\]
For example, $F$ has height $0$ at $p$ if and only if the field $k$
does not have characteristic $p$.  Over an algebraically closed field
of characteristic $p$, the height invariant $ht_p$ determines the
formal group up to isomorphism (but {\em not} up to strict
isomorphism).

\begin{thm}[\cite{landweber, rudyak}]
Suppose that $M$ is a graded module over the Lazard ring $L$.  Then the
functor sending a spectrum $X$ to the graded abelian group
\[
M \otimes_L MU_*(X)
\]
defines a generalized homology theory if and only if, for all primes
$p$ and all $n$, the map $v_n$ is an injective self-map of
$M/(p,\ldots,v_{n-1}).$
\end{thm}

We refer to such an object as a {\em Landweber exact} theory.  May
showed that such theories can be realized by $MU$-modules
\cite[Theorem 8]{may}, and Hovey-Strickland showed that there is a
functorial lifting from the category of Landweber exact theories to
the homotopy category of $MU$-modules \cite{hoveystrickland}.  In
addition, there are results for $L$-algebras rather than $L$-modules.

This theorem can be used to gives rise to numerous theories; complex
$K$-theory $KU$ is one such by the Conner-Floyd theorem.  Other examples
include the Brown-Peterson spectra $BP$ and Johnson-Wilson spectra $E(n)$.

In the case of complex $K$-theory, we have also have a more refined
multiplicative structure and the Adams operations $\psi^r$.  There is
a generalization of this structure due to Goerss-Hopkins-Miller
\cite{rezknotes,goersshopkins}.

Associated to a formal group law $F$ over a perfect field $k$ of
characteristic $p$, there is a complete local ring $\LT(k,F)$, called
the Lubin-Tate ring, with residue field $k$.  The Lubin-Tate ring
carries a formal group law $\tilde F$ equipped with an isomorphism of
its reduction with $F$.  If $F$ has $ht_p(F) = n$, then
\[
\LT(k,F) \cong \mb W(k)\llbracket u_1,\cdots, u_{n-1}\rrbracket,
\]
where $\mb W(k)$ is the Witt ring of $k$.

This ring is universal among such local rings, as follows.  Given any
local ring $R$ with nilpotent maximal ideal ${\mf m}$ and residue
field an extension $\ell$ of $k$, together with a formal group law $G$
over $R$ such that $G$ and $F$ have the same extension to $\ell$,
there exists a unique ring map $\LT(k,F) \to R$ carrying $\tilde F$ to
$G$.  In particular, the group of automorphisms of $F$ acts on
$\LT(k,F)$.

\begin{thm}[Goerss-Hopkins-Miller]
There is a functor
\[
E\co \{\text{formal groups over perfect fields, isos}\} \to \{E_\infty \text{ ring
  spectra}\}
\]
such that the homotopy groups of $E(k,F)$ are $\LT(k,F)[u^{\pm 1}]$,
where $|u| = 2$.
\end{thm}

This spectrum is variously referred to as a Hopkins-Miller spectrum,
Lubin-Tate spectrum, or Morava $E$-theory spectrum.  It is common to
denote by $E_n$ the spectrum associated to the particular example of
the Honda formal group law over the field $\mb F_{p^n}$, which has
height $n$.  Even worse, this theory is sometimes referred to as {\em
  the} Lubin-Tate theory of height $n$.  To do so brushes the
abundance of different multiplicative forms of this spectrum under the
rug.

We note that this functorial behavior allows us to construct
cohomology theories that are not complex oriented.  For instance, the
real $K$-theory spectrum $KO$ is the homotopy fixed point spectrum of
the action of the group $\{1,\psi^{-1}\}$ on $KU$, and the
$K(n)$-local spheres $L_{K(n)} \mb S$ are fixed point objects of
the full automorphism groups of the Lubin-Tate theories
\cite{devinatzhopkins}.

The extra multiplicative structure on the Lubin-Tate spectra allows us
to speak of categories of modules and smash products over them, both
powerful tools in theory and application.  The functoriality in the
Goerss-Hopkins-Miller theorem allows one to construct many new spectra
via homotopy fixed-point constructions.  These objects are now
indispensable in stable homotopy theory.

\section{Forms of the multiplicative group}
\label{sec:multiplicative}

The purpose of this section is to describe real $K$-theory as being
recovered from families of formal group laws, and specifically
cohomology theories associated to forms of the multiplicative group.

There is a multiplicative group scheme $\mb G_m$ over $\mb Z$.  It
is described by the Hopf algebra $\mb Z[x^{\pm 1}]$, with comultiplication
$x \mapsto x \otimes x$.  For a ring $R$, the set of $R$-points of
$\mb G_m$ is the unit group $R^\times$.  The formal completion of this
at $x=1$ is a formal group $\hat{\mb G}_m$.

However, there are various nonisomorphic {\em forms} of the
multiplicative group over other base rings that become isomorphic
after a flat extension.  For example, there is a Hopf algebra
\[
\mb Z\left[\frac{1}{2},x,y\right]/(x^2 + y^2 - 1),
\]
with comultiplication $x \mapsto (x \otimes x - y \otimes y), y
\mapsto (x \otimes y + y \otimes x)$.  For a ring $R$, the set of
$R$-points is the set
\[
\{x + iy\ |\ x^2 + y^2 = 1\},
\]
with multiplication determined by $i^2 = -1$.  Although all forms of
the multiplicative group scheme become isomorphic over an
algebraically closed field, there is still number-theoretic content
locked into these various forms.

We now parametrize these structures.  Associated to any pair of
distinct points $\alpha,\beta \in \mb A^1$, there is a unique group
structure on $\mb P^1 \setminus \{\alpha,\beta\}$ with $\infty$ as
unit.  The pair of points is determined uniquely by being the roots of
a polynomial $x^2 + bx + c$ with discriminant $\Delta = b^2 - 4c$
a unit.  Explicitly, the group structure is given by
\[
(x_1,x_2) \mapsto \frac{x_1 x_2 - c}{x_1 + x_2 - b}.
\]
This has a chosen coordinate $1/x$ near the identity of the group
structure.  By taking a power series expansion of the group law, we
get a formal group law.  We note that given $b$ and $c$ in a ring $R$,
we can explicitly compute the $p$-series as described in section
\ref{sec:realize}, and find that the image of $v_1 \in L/p$ is
\[
(\beta - \alpha)^{p-1} = \Delta^{\frac{p-1}{2}}.
\]
Therefore, such a formal group law over a ring $R$ is always Landweber
exact when multiplication by $p$ is injective for all $p$.

An isomorphism between two such forms of $\mb P^1$ must be given by an
automorphism of $\mb P^1$ preserving $\infty$, and hence a linear
translation $x \mapsto \lambda x + r$.  Expanding in terms of $1/x$,
such an isomorphism gives rise to a {\em strict} isomorphism if and
only if $\lambda = 1$.

We therefore consider the following three Hopf algebroids
parametrizing isomorphism classes of quadratics $x^2 + bx + c$, or
forms of the multiplicative group, in different ways.
\begin{eqnarray*}
  A &=& \mb Z[b,c, (b^2 - 4c)^{-1}]\\
  \Gamma_A &=& A[r]\\
  B &=& \mb Z[\alpha, \beta, (\alpha - \beta)^{-1}]\\
  \Gamma_B &=& B[r,s]/(s^2 + (\alpha-\beta)s)\\
  C &=& \mb Z[\alpha^{\pm 1}]\\
  \Gamma_C &=& C[s]/(s^2 + \alpha s)  
\end{eqnarray*}
These determine categories such that, for any ring $T$, the $T$-points
are given as follows.

\begin{tabular}{cl}
$(A,\Gamma_A)\co$ &  \{quadratics $x^2 + bx + c$, translations $x
\mapsto x+r$\}\\
$(B,\Gamma_B)\co$ & \{quadratics $(x-\alpha)(x-\beta)$,\\
&\ translations $x \mapsto x+r$ plus interchanges of $\alpha$ and $\beta$\}\\
$(C,\Gamma_C)\co$ & \{quadratics $x^2 - \alpha x$, transformations $x
\mapsto x + \alpha$\}\\
\end{tabular}

There is a natural faithfully flat map $A \to B$ given by $b \mapsto
-(\alpha + \beta), c \mapsto \alpha \beta$ corresponding to a
forgetful functor on quadratics.  The induced descent Hopf algebroid
$(B, B \otimes_A \Gamma_A \otimes_A B)$ is isomorphic to $(B,
\Gamma_B)$, and so the two Hopf algebroids represent the same stack.

The category given by the second is naturally equivalent to a
subcategory given by the third Hopf algebroid for all $T$.  We can choose
a universal representative for this natural equivalence given by the
natural transformation $\Gamma_B \to B$ of $B$-algebras sending
$r$ to $\beta$ and $s$ to $0$, showing that these also represent the
same stack.

The third Hopf algebroid, finally, is well-known as the Hopf algebroid
computing the homotopy of real $K$-theory $KO$.

In this way, we ``recover'' real $K$-theory as being associated to the
moduli stack of forms of the multiplicative group in a way compatible
with the formal group structure.

We note that, by {\em not} inverting the discriminant $b^2 -4c$, we
would recover a Hopf algebroid computing the homotopy of the {\em
  connective} real $K$-theory spectrum $ko$.  On the level of moduli
stacks, this allows the degenerate case of the additive formal group
scheme $\mb G_a$ of height $\infty$.  Geometrically, this point is
dense in the moduli of forms of $\mb G_m$.

\section{Elliptic curves and elliptic cohomology theories}
\label{sec:elliptic}

One other main source of formal group laws in algebraic geometry is
given by elliptic curves.

Over a ring $R$, any equation of the form
\[
y^2 + a_1 xy + a_3 y = x^3 + a_2 x^2 + a_4 x + a_6
\]
(a Weierstrass equation) determines a closed subset of projective
space $\mb P^2$.  There is a discriminant invariant $\Delta \in R$
which is a unit if and only if the group scheme is smooth.  See
\cite[Chapter III]{silverman}.

There is a commutative group law on the nonsingular points with
$[0:1:0] \in \mb P^2$ as identity.  Three distinct points $p$, $q$,
and $r$ are colinear in $\mb P^2$ if and only if they add to zero in
the group law.

The coordinate $x/y$ determines a coordinate near $\infty$ in the
group scheme, and expanding the group law in power series near
$\infty$ gives a formal group law over $R$.

Two Weierstrass curves are isomorphic over $R$ if and only if there is
a unit $\Lambda \in R^\times$ and $r,s,t \in R$ such that the
isomorphism is given by $x \mapsto \lambda^2 x + r, y \mapsto
\lambda^3 y + sx + t$.  The isomorphism induces a strict isomorphism
of formal group laws if and only if $\lambda = 1$.

An elliptic curve over a general scheme has a formal definition, but
can be formed by patching together such Weierstrass curves locally (in the
flat topology).  There is a Hopf algebroid representing the
groupoid of nonsingular Weierstrass curves and strict isomorphisms, given by
\begin{eqnarray*}
  A &=& \mb Z[a_1,a_2,a_3,a_4,a_6, \Delta^{-1}],\\
  \Gamma &=& A[r,s,t].
\end{eqnarray*}
The associated stack ${\cal M}_{ell}$ is a moduli stack of elliptic
curves (and strict isomorphisms).  The natural association taking such
an elliptic curve to its formal group law gives a map of stacks
\[
{\cal M}_{ell} \to {\cal M}_{sFG}
\]
to the moduli stack of formal groups.

One can instead think of this moduli stack as parametrizing pairs
$(E,\omega)$ of an elliptic curve $E$ and a nonzero invariant $1$-form
$\omega$ on $E$.  The invariant $1$-form determines a coordinate near
the unit of the elliptic curve up to first order, and a map of such
elliptic curves then induces a strict isomorphism if and only if it
preserves the form.

An {\em elliptic cohomology theory} consists of a cohomology theory
$E$ which is weakly even periodic\footnote{A spectrum is {\em weakly
    even periodic} if the nonzero homotopy groups are concentrated in
  even degrees, and the product $E_p \otimes_{E_0} E_q \to E_{p+q}$ is
  always an isomorphism for $p,q$ even.}, together with an elliptic
curve over $\Spec(E_0)$ and an isomorphism of formal group laws
between the formal group law associated to the elliptic curve and the
formal group law of the spectrum.  Landweber exact theories of this
form were investigated by Landweber-Ravenel-Stong based on a Jacobi
quartic \cite{lrs}.  In terms of the moduli, we would like to view
these as arising from schemes $\Spec(E_0)$ over ${\cal M}_{ell}$ with
spectra realizing them.

Similarly, by allowing the possibility of elliptic curves with {\em
nodal} singularities (so that the resulting curve is isomorphic to
$\mb P^1$ with two points identified, with multiplication on the
smooth locus a form of $\mb G_m$), we get a compactification
$\overline {{\cal M}_{ell}}$ of the moduli of elliptic curves.  The
object $\overline {{\cal M}_{ell}}$ is a smooth Deligne-Mumford stack
over $\Spec(\mb Z)$ \cite{delignerapoport}.  This stack is more
difficult to express in terms of Hopf algebroids.

Based on our investigation of forms of the multiplicative formal
group, it is natural to ask whether there is a ``universal'' elliptic
cohomology theory associated to ${\cal M}_{ell}$ and a universal
elliptic cohomology theory with nodal singularities associated to
$\overline {{\cal M}_{ell}}$.  Here we could interpret universality as
being either a lift of the universal elliptic curve over this stack,
or being somehow universal among elliptic cohomology theories.

If $6$ is invertible in $R$, each Weierstrass curve is isomorphic (via
a unique {\em strict} isomorphism) to a uniquely determined elliptic curve
of the form $y^2 = x^3 + c_4 x + c_6$.  This universal elliptic curve
over the (graded) ring $\mb Z[\frac{1}{6},c_4,c_6,\Delta^{-1}]$ has a
Landweber exact formal group law, and hence is realized by a
cohomology theory generally denoted by $\Ell$ \cite{baker}.

We would be remiss if we did not mention the inspiring connection to
multiplicative genera and string theory \cite{ahs}.

\section{Topological modular forms}
\label{sec:tmf}

The theories $\TMF[\Delta^{-1}]$, $\TMF$, and $\tmf$ of topological
modular forms are extensions of the construction of the universal
elliptic theory $\Ell$.  This extension occurs in several directions.
\begin{itemize}
\item These theories are all realized by $E_\infty$ ring spectra, with
  the corresponding increase in structure on categories of modules and
  algebras.
\item These theories are universal objects, in that they can be
  constructed as a limit of elliptic cohomology theories.
  $\TMF[\Delta^{-1}]$ and $\TMF$ are associated to the moduli stacks
  ${\cal M}_{ell}$ and $\overline {{\cal M}_{ell}}$ respectively.  These
  are not elliptic cohomology theories themselves, just as $KO$ is not
  a complex oriented theory due to the existence of forms of $\mb G_m$
  with automorphisms.
\item Unlike $\Ell$, these theories carry information at the primes
  $2$ and $3$.  In particular, they detect a good portion of
  interesting $2$- and $3$-primary information about stable homotopy
  groups of spheres.
\end{itemize}
The construction of these theories (due to Hopkins et al.) has yet to
fully appear in the literature, but has nevertheless been highly
influential in the subject for several years.

An interpretation in terms of sheaves is as follows.  On the moduli
${\cal M}_{ell}$ and $\overline {{\cal M}_{ell}}$ of elliptic curves,
any \'etale map (roughly, a map which is locally an isomorphism, such
as a covering map) from $\Spec(R)$ can be realized by a highly
structured elliptic cohomology theory in a functorial way.  Stated
another way, we have a lift of the structure sheaf ${\cal O}$ of the
stack in the \'etale topology to a sheaf ${\cal O}^{der}$ of
commutative ring spectra.

(We should mention that associated to {\em modular curves}, which are
certain coverings of $\overline {{\cal M}_{ell}}$, these structure sheaves give
rise to versions of $\TMF$ with level structures.  This construction,
however, may require certain primes to be inverted.)

The homotopy of $\TMF[\Delta^{-1}]$ is computable via the
Adams-Novikov spectral sequence \cite{bauer, rezkcourse}, whose
$E_2$-term is the cohomology of the Weierstrass curve Hopf algebroid
of section \ref{sec:elliptic}.  Similarly, the Adams-Novikov spectral
sequence for the homotopy of $\TMF$ has $E_2$-term given in terms of
the cohomology of the compactified moduli $\overline {{\cal
    M}_{ell}}$.  The zero-line of each of these spectral sequences can
be identified with a ring of modular forms over $\mb Z$.

The spectrum $\tmf$ also has homotopy computed by the Weierstrass
algebroid, but without the discriminant inverted.  It corresponds to a
moduli of possibly singular elliptic curves where we allow the
possibility of curves with {\em additive} reduction, or cusp
singularities.  As a spectrum, however, $\tmf$ is generally
constructed as a connective cover of $\TMF$ and does not fit well into
the theory of ``derived algebraic geometry'' due to Lurie et al.

\section{The moduli stack of formal groups}
\label{sec:chromatic}

We have discussed several cohomology theories here with relationships
to the moduli stack of formal groups ${\cal M}_{sFG}$.  It is time to
elaborate on the geometry of this moduli stack.

{\em From this point forward, we fix a prime $p$ and focus our
  attention there.  In particular, all rings and spectra are assumed
  to be $p$-local, or $p$-localized if not.}

We recall that a formal group law over an algebraically closed field
of characteristic $p$ is classified uniquely up to isomorphism by its
height invariant.  In terms of the Lazard ring $L$, we have a sequence
of elements $p,v_1,v_2,\cdots$, with each prime ideal
$(p,v_1,\ldots,v_{n-1})$ cutting out an irreducible closed substack ${\cal
  M}_{sFG}^{\geq n}$ of the moduli stack.  It turns out that these
prime ideals (and their union) are the {\em only} invariant prime
ideals of the moduli.  The intersection of all these closed substacks
is the height-$\infty$ locus.

As a result, we have a stratification of the moduli stack into layers
according to height.  There is a corresponding filtration in homotopy
theory called the {\em chromatic filtration}, and it has proved to be
a powerful organizing principle for understanding large-scale
phenomena in homotopy theory \cite{ravenelgreen,dhs}.  We note
that the Landweber exact functor theorem might be interpreted as a
condition for a map $\Spec(R) \to {\cal M}_{sFG}$ to be flat.

Having said this, we would like to indicate how the various cohomology
theories we have discussed fit into this filtration.

Rational cohomology, represented by the Eilenberg-Maclane spectrum
$\eilm{\mb Q}$, has the prime $p = v_0$ inverted.  It hence lives over
the height $0$ open substack of ${\cal M}_{sFG}$.

Mod-$p$ cohomology, represented by the Eilenberg-Maclane spectrum
$\eilm{\mb F_p}$, has the additive formal group law $x +_F y = x + y$,
and hence is concentrated over the height $\infty$ closed substack.

We saw in section \ref{sec:multiplicative} that forms of the
multiplicative formal group law have the quantity $v_1$ invertible.
These theories, exemplified by complex $K$-theory $KU$ and real
$K$-theory $KO$, therefore are concentrated over the open substack of
heights less than or equal to $1$.  (The connective versions $ko$ and
$ku$ of these spectra are concentrated over heights $0$, $1$, and
$\infty$.)  The work of Morava on forms of $K$-theory also falls into
this region \cite{morava}.

It is a standard part of the theory of elliptic curves in
characteristic $p$ that there are two distinct classes: the {\em
  ordinary} curves, whose formal groups have height 1, and the {\em
  supersingular} curves, whose formal groups have height 2.  The
theories $\TMF[\Delta^{-1}]$ and $\TMF$, and indeed all elliptic
cohomology theories, are therefore concentrated on the open substack
of heights less than or equal to $2$.  (The connective spectrum $\tmf$
is concentrated over heights $0$, $1$, $2$, and $\infty$.)

As these theories only detect ``low'' chromatic phenomena, they are
limited in their ability to detect phenomena in stable homotopy
theory.  It is natural to ask for us to find cohomology theories that
elaborate on the chromatic layers in homotopy theory at all heights.

It is worth remarking that an understanding of chromatic level one led
to proofs of the Hopf invariant one problem, and hence to the final
solution of the classical problem about vector fields on spheres.
Referring to chromatic level two as ``low'' is incredibly misleading.
The computations involved in stable homotopy theory at chromatic level
two are quite detailed \cite{shimomurawang,ghmr}, and the Kervaire
invariant problem is concentrated at this level. Very little is
computationally known beyond this point.

Several examples of spectra with higher height are given by the
Morava theories mentioned in section \ref{sec:realize}.  The Morava
$E$-theory spectrum $E(k,F)$ associated to a formal group law of
height $n < \infty$ over a perfect field $k$ is concentrated over the
height $\leq n$ open substack of ${\cal M}_{sFG}$.  In some sense,
however, these theories are controlled by their behavior at height
exactly $n$, and do not have much ``interpolating'' behavior.  They
are also more properly viewed as ``pro-objects'' (inverse systems) in
the stable homotopy category, and have homotopy groups that are not
finitely generated as abelian groups.  Finally, these theories are
derived strictly from the formal group point of view in homotopy
theory, and they can be difficult to connect to geometric content.

More examples are given by the Johnson-Wilson theories $E(n)$, which
are not known to have much structured multiplication for $n > 1$.

More ``global'' examples are given by spectra denoted $\eo_n$, where
$\eo_2$ is $\tmf$.  These spectra take as starting point the
Artin-Schreier curve
\[
y^{p-1} = x^p - x.
\]
In characteristic $p$, this curve has a large symmetry group that also
acts on the Jacobian variety.  The Jacobian has a higher-dimensional
formal group, but the group action produces a 1-dimensional split
summand of this formal group with height $p-1$.  Hopkins and
Gorbunov-Mahowald\footnote{The author's talk at the conference
misattributed this, and multiple attendees corrected him; he would
like to issue an apology.} initiated an investigation of a
Hopf algebroid associated to deformations of this curve of the form
\[
y^{p-1} = x^p - x + \sum u_i x^i,
\]
whose realization would be a spectrum denoted by $\eo_{p-1}$
\cite{gorbunovmahowald}.  Ravenel generalized this to the
Artin-Schreier curve
\[
y^{p^f - 1} = x^p - x,
\]
whose formal group law has a 1-dimensional summand of height $(p-1)f$
and an interesting symmetry group \cite{ravenelhigher}.  However, the
existence of spectrum realizations is (at the time of this writing)
still not known.

\section{$p$-divisible groups and Lurie's theorem}
\label{sec:lurie}

In 2005, Lurie announced a result that gave sufficient conditions to
functorially realize a family of 1-dimensional formal group laws by
spectra given certain properties and certain extra data.  The extra
data comes in the form of a $p$-divisible group (or Barsotti-Tate
group), and the necessary property is that locally the structure of
the $p$-divisible group determines the geometry.  In this section we
introduce some basics on these objects.  The interested reader should
consult \cite{tateflat} or \cite{messing}.

A {\em $p$-divisible group} $\mb G$ over a an algebraically closed
field $k$ consists of a (possibly multi-dimensional) formal group $F$
of finite height $h$ and a discrete group isomorphic to $(\mb Q/\mb
Z)^r$, together in an exact sequence
\[
0 \to F \to \mb G \to (\mb Q/\mb Z)^r \to 0.
\]
The integer $n = h+r$ is the {\em height} of $\mb G$, and the
dimension of the formal component $F$ is the {\em dimension} 
of $\mb G$.

However, we require a more precise description in general.  Over a
base scheme $X$, a $p$-divisible group actually consists of a sequence
of finite, flat group schemes $\mb G[p^k]$ (the $p^k$-torsion) over
$X$ with $\mb G[p^0] = 0$ and inclusions $\mb G[p^k] \subset \mb
G[p^{k+1}]$ such that the multiplication-by-$p$ map factors as
\[
\mb G[p^{k+1}] \onto \mb G[p^k] \subset \mb G[p^{k+1}].
\]
The height and dimension of the $p$-divisible group are locally
constant functions on $X$, equivalent to the rank of $\mb G[p]$ and
the dimension of its tangent space.  At any geometric point $x \in X$,
the restriction of the $p$-divisible group to $x$ lives in the natural
short exact sequence
\[
0 \to \mb G^{for} \to \mb G_x \to \mb G^{\text{\'et}} \to 0,
\]
with the subobject (the connected component of the unit) the
{\em formal} component and the quotient the {\em \'etale} component.  The
formal component $\mb G^{for}$ is a formal group on $X$.  The height
of the formal component is an upper semicontinuous function on $X$,
and gives rise to a stratification of $X$ which is the pullback of the
stratification determined by the regular sequence $(p,v_1,\ldots)$.

In fact, a deeper investigation into the isomorphism classes of
$p$-divisible groups over a field gives rise to a so-called ``Newton
polygon'' associated to a $p$-divisible group and a Newton polygon
stratification.  However, for $p$-divisible groups of dimension $1$
this is equivalent to the formal-height stratification.

Similar to formal group laws, there is a deformation theory of
$p$-divisible groups.  Each $p$-divisible group $\mb G$ of height $n$
over a perfect field $k$ of characteristic $p$ has a universal
deformation $\tilde {\mb G}$ over a ring analogous to the Lubin-Tate
ring.

For any $n < \infty$, there is a formal moduli ${\cal M}_p(n)$ of
$p$-divisible groups of height $n$ and dimension $1$ and their
isomorphisms.  The author is not aware of any amenable presentations
of a moduli stack analogous to the presentation of the moduli of
formal group laws, and whether a well-behaved Hopf algebroid exists
modelling this stack seems to still be open.  From a formal point of
view, the category of maps from a scheme $X$ to ${\cal M}_p(n)$ should
be the category of $p$-divisible groups of height $n$ on $X$, and the
association $\mb G \mapsto \mb G^{for}$ gives a natural transformation
from ${\cal M}_p(n)$ to the moduli of formal groups ${\cal M}_{FG}$.

We state a version of Lurie's theorem here.

\begin{thm}[Lurie]
  Let ${\cal M}$ be an algebraic stack over $\mb Z_p$\footnote{The
  stack ${\cal M}$ must actually be formal, with $p$ topologically
  nilpotent.} equipped with a morphism
\[
{\cal M} \to {\cal M}_p(n)
\]
classifying a $p$-divisible group $\mb G$.  Suppose that at any point
$x \in {\cal M}$, the complete local ring of ${\cal M}$ at $x$ is
isomorphic to the universal deformation ring of the $p$-divisible
group at $x$.  Then the composite realization problem
\[
{\cal M} \to {\cal M}_p(n) \to {\cal M}_{FG}
\]
has a canonical solution; that is, there is a sheaf of $E_\infty$ even
weakly periodic $E$ on the etale site of ${\cal M}$ with $E_0$ locally
isomorphic to the structure sheaf and the associated formal group
${\cal G}$ isomorphic to the formal group $\mb G^{for}$.  The space of
all solutions is connected and has a preferred basepoint.
\end{thm}

The proof of Lurie's theorem requires the Hopkins-Miller theorem to
provide objects for local comparison, and so generalizations without
the ``universal deformation'' condition are not expected without some
new direction of proof.  We also note that the theorem does not apply
as stated to the compactified moduli $\overline{{\cal M}_{ell}}$, and
so only gives a proof of the existence of $\TMF[\Delta^{-1}]$ rather
than $\TMF$.

Our perspective, however, is to view this theorem as a black box.  It
tells us that if we can find a moduli ${\cal M}$ such that
\begin{itemize}
\item ${\cal M}$ has a {\em canonically} associated 1-dimensional
  $p$-divisible group $\mb G$ of height $n$, and
\item the local geometry of ${\cal M}$ corresponds exactly to local
  deformations of $\mb G$,
\end{itemize}
then we can find a canonical sheaf of spectra on ${\cal M}$.  Having
this in hand, our goal is to seek examples of such moduli.

Unfortunately, several examples mentioned in previous sections do not
immediately seem to have attached $p$-divisible groups.  The
deformations of Artin-Schreier curves in the previous section, or
Johnson-Wilson theories, do not {\em a priori} have attached
$p$-divisible groups.\footnote{In the Artin-Schreier case, the
  question becomes one of deforming the 1-dimensional split summand of
the Jacobian at the Artin-Schreier curve to a 1-dimensional
$p$-divisible group at all points.  The author is not aware of a
solution to this problem at this stage.}

At the other extreme, one could ask to realize the moduli stack ${\cal
  M}_p(n)$ itself by a spectrum.  This stack has geometry very close
to the moduli of formal groups.  In particular, it still breaks down
according to height, but is truncated at height $n$ and has extra
structure at heights below $n$.  The resulting object should give an
interesting perspective on chromatic homotopy theory.

The main obstruction to this program, however, seems to be the
difficulty in finding a presentation of this stack or any reasonable
information about the category of quasicoherent sheaves.

\section{PEL Shimura varieties and $\TAF$}
\label{sec:pel}

Based on Lurie's theorem, it becomes natural to seek moduli problems
with associated $1$-dimensional $p$-divisible groups of height $n$ in
order to produce new spectra.  Following the approaches of
Gorbunov-Mahowald and Ravenel, we approach this through abelian
varieties.  However, rather than considering families of plane curves
and their Jacobians, we consider families of abelian varieties
equipped with extra structure.  The stunning fact is that the precise
assumptions needed to produce reasonable families of $p$-divisible
groups occur {\em already} in families of PEL abelian varieties of a
type studied classically by Shimura, and of the specific kind featured
in Harris and Taylor's proof of the local Langlands correspondence
\cite{harristaylor}.  The reader interested in these varieties should
refer to \cite{milne} and then \cite{kottwitz}.

One of the main places that $p$-divisible groups occur in algebraic
geometry is from group schemes.  For any (connected) commutative group
scheme $G$, we have maps representing multiplication by $p^k$:
\[
[p^k]\co G \to G.
\]
The identity element $e \in G$ has a {\em scheme-theoretic} inverse
image $G[p^k] \subset G$.  Associated to a group scheme $G$ over a
given base $X$, the system $G[p^k]$ forms a $p$-divisible group $G(p)$
under sufficient assumptions on $G$, such as if $G$ is an abelian
variety.

For example, consider the multiplicative group scheme over $\Spec(R)$,
given by $\mb G_m = \Spec(R[t^{\pm 1}])$.  The multiplication-by-$p^k$ map is given on
the ring level by the map $t \mapsto t^{p^k}$, and the scheme-theoretic
preimage of the identity is the subscheme of solutions of $t^{p^k} =
1$, or
\[
\Spec(R[t^{\pm 1}]/(t^{p^k} - 1)).
\]
If $R$ has characteristic zero, then this scheme has $p^k$ distinct
points over each geometric point of $\Spec(R)$.  If $R$ has
characteristic $p$, then this scheme is isomorphic to
\[
\Spec(R[t^{\pm 1}]/(t-1)^{p^k}).
\]
Each geometric point has only {\em one} preimage in this case, and so
the $p$-divisible group $\mb G_m(p)$ is totally formal.

The basic problem is as follows.
\begin{itemize}
\item The only $1$-dimensional group schemes over an algebraically
  closed field are the additive group $\mb G_a$, the multiplicative
  group $\mb G_m$, and elliptic curves.
\item The $p$-divisible group of an $n$-dimensional abelian variety
  $A$ has height $2n$ and dimension $n$.
\end{itemize}
As a result, if we decide that we will consider moduli of
higher-dimensional abelian varieties, we need some way to cut down the
dimension of the $p$-divisible group to 1.  As in the
Mahowald-Gorbunov-Ravenel approach, we can carry this out by assuming
that we have endomorphisms of the abelian variety splitting off a
1-dimensional summand $\mb G$ canonically.

However, we also must satisfy a condition on the local geometry.  What
this translates to in practice is the following: given an infinitesimal
extension of the $p$-divisible group $\mb G$, we must be able to
complete this to a unique deformation of the element in the moduli.

Our main weapon in this task is the following.  See \cite{katz} for a
proof, due to Drinfel'd.

\begin{thm}[Serre-Tate]
Suppose we have a base scheme $X$ in which $p$ is locally nilpotent,
together with an abelian scheme\footnote{An abelian scheme is a
family of abelian varieties over the base.} $A/X$.  Any
deformation of the $p$-divisible group $A(p)$ determines a unique
deformation of $A$.
\end{thm}

Some of the language here is deliberately vague.  However, this is
more easily stated in terms of fields.  Suppose that $k$ is a field of
characteristic $p$, and $R$ is a local ring with nilpotent maximal
ideal ${\mf m}$ and residue field $k$.  Then the category of abelian
schemes over $R$ is naturally equivalent (via a forgetful functor) to
the category of abelian varieties $A$ over $k$ equipped with
extensions of their $p$-divisible group $A(p)$ to $R$.

This does the heavy lifting for us.  If we can specify a moduli of
abelian varieties with a $1$-dimensional summand $\mb G$ of the
$p$-divisible group that controls the {\em entire} $p$-divisible group
in some way, we will be done.  This is accomplished via the
aforementioned moduli of PEL Shimura varieties.  For simplicity, we
consider the case of simple complex multiplication, rather than an
action by a division algebra, leaving generality to other references.

To define these Shimura varieties requires the compilation of a
substantial dossier.  We simply present this now, and make it our goal
in the following sections to justify why all these pieces of data are
important for us to include.

We first must state some necessary facts from the theory of abelian
schemes without proof.
\begin{itemize}
\item If $A$ is an abelian scheme, the {\em dual} abelian scheme
  $A^\vee$ is the identity component ${\rm Pic}^0(A)$ of the group of
  line bundles on $A$.  Duals exist over a general base scheme,
  dualization is a contravariant functor, and the double-dual is
  canonically isomorphic to $A$.
\item There is a compatible dualization functor on $p$-divisible
  groups with a canonical isomorphism $A^\vee(p) \cong (A(p))^\vee$.
  Dualization preserves height, but not dimension.  However, we have
  that $\dim(\mb G) + \dim(\mb G^\vee)$ is the height of $\mb G)$.
\item An isogeny $A \to B$ between abelian schemes is a surjection
  with finite kernel; it expresses $B$ as isomorphic to $A/H$ for $H$
  a finite subgroup scheme of $A$.  An isogeny is {\em prime-to-$p$}
  if the kernel has rank prime to $p$ (as a group scheme).
\item The endomorphism ring $\End(A(p))$ is $p$-complete, and hence a
  $\mb Z_p$-algebra.
\end{itemize}

Fix an integer $n$ and continue to fix a prime $p$.  Let $F$ be a
quadratic imaginary extension field of $\mb Q$, and ${\cal O}_F$ the
ring of integers of $F$.  We require that $F$ be chosen so that 
$p$ splits in $F$, i.e. ${\cal O}_F \otimes \mb Z_p \cong \mb Z_p
\times \mb Z_p$.  In particular, we can choose an idempotent $e \in
{\cal O}_F \otimes \mb Z_p$ such that $e \neq 0,1$.  Complex
conjugation is forced to take $e$ to $1-e$.

In addition, we need to fix one further piece of data required to
specify a level structure, which will be discussed in section
\ref{sec:level}.

We consider the functor that associates to a scheme $X$ over
$\mb Z_p$ the category of tuples $(A,\lambda, \iota, \eta)$ of the
following type.
\begin{itemize}
\item $A$ is an abelian scheme of dimension $n$ over $X$.
\item $\lambda\co A \to A^\vee$ is a prime-to-$p$ polarization.  (This
  is an isogeny such that $\lambda^\vee = \lambda$, together with a
  positivity condition; we will discuss it further in section
  \ref{sec:pol}.)
\item $\iota\co {\cal O}_F \to \End(A)$ is a ring homomorphism from
  ${\cal O}_F$ to the endomorphism ring of $A$ such that $\lambda
  \iota(\alpha) = \iota(\bar \alpha)^\vee \lambda$ for all $\alpha \in
  {\cal O}_F$.  We require that the summand $e \cdot A(p) \subset
  A(p)$ is $1$-dimensional.  (See section \ref{sec:end}.)
\item $\eta$ is a level structure on $A$.  (See section
  \ref{sec:level}.)
\end{itemize}

Morphisms in the category are isomorphisms $f\co A \to B$ that
commute with the action $\iota$, that preserve the level structure,
and such that $f^\vee \lambda_B f = n \lambda_A$ for some positive
integer $n$.

We take as given that this moduli is well-behaved.  In particular, it
is represented by a smooth Deligne-Mumford stack of relative dimension
$(n-1)$ over $\mb Z_p$.  We abusively denote it by $\Sh$ without
decorating it with {\em any} of the necessary input data.  It has an
associated sheaf of spectra, and the ``universal'' object (a limit, or
global section object) is denoted $\TAF$.  The Adams-Novikov spectral
sequence takes the form
\[
H^s(\Sh, \omega^{\otimes t}) \Rightarrow \pi_{t-s} \TAF,
\]
where $\omega$ is the line bundle of invariant $1$-forms on the
1-dimensional formal component.  The zero line
\[
H^0(\Sh, \omega^{\otimes t})
\]
consists of (integral) automorphic forms on the Shimura stack.

The height $n$ stratum of the Shimura stack is nonempty, and consists
of a finite set of points whose automorphism groups can be identified
with finite subgroups of the so-called Morava stabilizer group $\mb
S_n$.  There is a corresponding description of the $K(n)$-localization
of the spectrum $\TAF$ as a finite product of fixed-point spectra of
Morava $E$-theories by finite subgroups.  These points can be
classified via the Tate-Honda classification of abelian varieties
over finite fields.

In the following sections, we will explain how the specified list of
data produces a 1-dimensional $p$-divisible group of the type
precisely necessary for Lurie's theorem.  For reasons of clarity in
exposition, we will discuss endomorphisms before polarizations.

\section{E is for Endomorphism}
\label{sec:end}

The most immediately relevant portion of the data of a Shimura variety
is the endomorphism structure.  The goal of this endomorphism is to
provide us with a 1-dimensional split summand of the $p$-divisible
group of $A$.

Recall that the endomorphism structure $\iota$ is a ring map ${\cal
O}_F \to \End(A)$, where ${\cal O}_F$ was a ring of integers whose
$p$-completion ${\cal O}_F \otimes \mb Z_p$ contains a chosen
idempotent $e$ making it isomorphic to $\mb Z_p \times \mb Z_p$.

The composite ring homomorphism
\[
{\cal O}_F \to \End(A) \to \End(A(p))
\]
lands in a $\mb Z_p$-algebra, and so we have a factorization
\[
{\cal O}_F \otimes \mb Z_p \to \End(A(p)).
\]
The image of the idempotent $e$ gives a splitting of $p$-divisible
groups
\[
A(p) \cong e\cdot A(p) \oplus (1-e) \cdot A(p).
\]
By assumption the $p$-divisible group $e \cdot A(p)$ is
1-dimensional.

Therefore, the elements of this moduli have canonically associated
1-dimensional $p$-divisible groups.  We do not yet know that these
have height $n$.

There is a similar decomposition of the $p$-divisible group of the
dual abelian variety.  
\[
A^\vee(p) \cong e^\vee \cdot A(p)^\vee \oplus (1-e^\vee) \cdot A(p)^\vee.
\]

\section{P is for Polarization}
\label{sec:pol}

The next piece of necessary data is the prime-to-$p$ polarization
$\lambda\co A \to A^\vee$.  Although polarizations are typically used
in algebraic geometry to guarantee representability of various moduli
problems (and this is a side effect necessary for us, as well), in our
case the polarization also gives control over the complementary
summand of the $p$-divisible group.

The condition that this map is a prime-to-$p$ isogeny implies that the
induced map of $p$-divisible groups $\lambda\co A(p) \to A^\vee(p)$ is
an isomorphism.

The condition that $\lambda$ conjugate-commutes with the action of
${\cal O}_F$ in particular implies
\[
\lambda e = (1-e^\vee) \lambda.
\]
As a result, the isomorphism
\[
A(p) \longoverto^\sim A^\vee(p)
\]
decomposes into the pair of isomorphisms.
\[
e \cdot A(p) \longoverto^\sim (1-e^\vee) \cdot A^\vee(p) = ((1-e) \cdot
A(p))^\vee
\]
\[
(1-e) \cdot A(p) \longoverto^\sim e^\vee \cdot A^\vee(p) = (e \cdot A(p))^\vee
\]
As a result, the polarization provides us with a {\em canonical}
identification of $(1-e) \cdot A(p)$, the $(n-1)$-dimensional
complementary summand of the $p$-divisible group, with the object $(e
\cdot A(p))^\vee$, the {\em dual} of the $1$-dimensional summand of
interest to us.  As the summands corresponding to $e$ and $(1-e)$ must
then have the same height, the height of each individual factor is
$n$.

This allows us to check that the conditions of Lurie's theorem hold.
As stated in section \ref{sec:lurie}, we must check that an
infinitesimal extension of the $1$-dimensional $p$-divisible group $e
\cdot A(p)$ determines a unique extension of $A$, with
endomorphisms, and with polarization.

In brief, we sketch the necessary reasoning.
\begin{itemize}
\item An extension of $e \cdot A(p)$ determines a dual extension of
  $(e \cdot A(p))^\vee \cong (1-e) \cdot A(p)$.
\item Therefore, we have an extension of the whole $p$-divisible
  group $A(p)$.
\item Declaring that $e$ and $(1-e)$ are idempotents corresponding to
  this splitting determines an extension of the action of ${\cal O}_F$.
\item The isomorphisms given by the polarization give a unique
  extension of $\lambda\co A(p) \to A(p)$ which conjugate-commutes
  with the action of ${\cal O}_F$.
\item The Serre-Tate theorem discussed in section \ref{sec:pel} then
  implies that the extension of $A(p)$, with the given extensions of
  $\iota$ and $\lambda$, determine a unique extension of $A$ with
  extensions of $\iota$ and $\lambda$.
\end{itemize}

A polarization also includes a positivity condition.  For a complex
torus $\mb C^g/\Lambda$ over $\mb C$, this amounts to a positive
definite Hermitian form on $\mb C^g$ whose imaginary part takes
integer values on $\Lambda$.  The existence of such a form serves to
eliminate the possibility that the torus does not have enough
nonconstant meromorphic functions on it to determine a projective
embedding; in higher dimensions, complex torii generically cannot be
made algebraic.

Polarizations also serve to eliminate pathology in families of abelian
varieties.  The set of automorphisms of a polarized abelian variety is
a finite group, and the moduli of polarized abelian varieties is
itself a Deligne-Mumford stack \cite{fc,mumford}.  Knowing this serves
as a first step in our ability to find a Deligne-Mumford stack for the
PEL moduli we are interested in.

\section{L is for Level Structure}
\label{sec:level}

There is one remaining ingredient in the data of a PEL Shimura
variety, which is the data of a level structure.

Those familiar with the more classical theory of elliptic curves will
be familiar with level structures such as the choice of a finite
subgroup of the curve, or a basis for the $n$-torsion.  This kind of
data can be included in the level structure, but it is not (for the
purposes of this document) the main point.

Given just the requirements of a polarization and endomorphism data (a
PE moduli problem), we would still have a moduli satisfying the
requirements of Lurie's theorem, and could produce spectra.  However,
such a moduli problem would usually suffer from a slight defect, in
the form of an infinite number of connected components.

There are various pieces of data, however, that are invariants of the
connected component; we can use this to classify various connected
components into ones of more manageable size for our sanity.

We require a definition.  Suppose $A$ is an abelian variety over an
algebraically closed field $k$.  For any prime $\ell \neq p$, we have
the groups $A[\ell^k]$ of $\ell$-torsion points of $A$, which are
abstractly isomorphic to $(\mb Z/\ell^k)^{2n}$.  These fit into an
inverse system
\[
\cdots \to A[\ell^3] \to A[\ell^2] \to A[\ell] \to 0
\]
where the maps are multiplication by $\ell$.  The inverse limit is
called the {\em $\ell$-adic Tate module} $T_\ell(A)$ of $A$, and is a
free $\mb Z_\ell$-module of rank $2n$.

The data of a polarization $A \to A^\vee$ gives rise to a pairing on
the $\ell$-adic Tate module.  Specifically, it gives rise to an
alternating bilinear pairing to the Tate module of the multiplicative
group scheme $T_\ell(\mb G_m) \cong \mb Z_\ell$.  This pairing is
referred to as the $\lambda$-Weil pairing.

If $(A,\lambda, \iota)$ is a polarized abelian variety over $k$ with
conjugate-commuting action of ${\cal O}_F$, we find that $T_\ell(A)$
is a free $\mb Z_\ell$-module of rank $2n$ equipped with a pairing
$\form{-}{-}$ on $\mb T_\ell(A)$.  This form is alternating, bilinear,
and ${\cal O}_F$-Hermitian in the sense that
\[
\form{\alpha x}{y} = \form{x}{\bar \alpha y}
\]
for all $\alpha \in {\cal O}_F$.

The isomorphism class of this pairing up to multiplication by a scalar
is an {\em invariant} of the connected component of $(A,\lambda, \iota)$ in
the PE moduli problem.

Therefore, part of the input data required to define our PEL moduli
problem is, for each $\ell \neq p$, a specified isomorphism class of
free $\mb Z_\ell$-module $M_\ell$ of rank $2n$ with alternating Hermitian
bilinear pairing (up to scale).  We can also specify an open subgroup
$K$ of the group of automorphisms of $\prod M_\ell$ (such as
automorphisms preserving specified subgroups or torsion points) as
part of the data.  The $K$-orbit of an isomorphism $\prod M_\ell \to
\prod T_\ell(A)$ is a {\em level $K$ structure}.

In the PEL moduli problem of tuples $(A,\lambda, \iota, \eta)$, the
level structure $\eta$ is a (locally constant) choice of level $K$
structure on $T_\ell(A_x)$ for each geometric point $x$ of the base
scheme $X$.  This is equivalent to specifying one such choice per
connected component which is invariant under the action of the \'etale
fundamental group of $X$.

Given such a level structure, one can prove that the moduli $\Sh$ over
$\mb Z_p$ consists of a finite number of connected components.  These
details do not occur in the elliptic case because there are few
isomorphism classes of alternating bilinear pairings on a lattice of
rank two.

It is common in the more advanced theory of automorphic forms to
simply drop the abelian varieties entirely, and simply think in terms
of a reductive algebraic group with a chosen open compact subgroup
$K$.  When pressed, for many expressions of a Shimura variety one can
find a reduction to a certain kind of moduli of abelian varieties by a
process of reduction.  However, this is by no means a straightforward
process.

\section{Questions}
\label{sec:questions}

This section is an attempt to give a series of straw-man arguments as
to why we might choose this particular conglomeration of initial data,
rather than making some slight alteration.  It also attempts to answer
some other questions that appear frequently.

\begin{ques}
Why do we act by ${\cal O}_F$ for a quadratic extension of $\mb Q$?
Why don't we choose endomorphisms by some other ring?  Why is $F$
specified as part of the data?
\end{ques}

In short, we must act by a ring whose $p$-completion contains an
idempotent, but does not contain an idempotent itself (which would
force the 1-dimensional summand to come from an elliptic curve, and
hence cap the height of the $p$-divisible group at 2).  In order to
uniquely give extensions of endomorphisms as in section \ref{sec:pol},
the $p$-completion of the ring must essentially be $\mb Z_p \times \mb
Z_p$, and since $\End(A)$ is a finitely generated free abelian group
for $A$ over a field $k$, we might as well assume that our ring to be
free of rank 2 over $\mb Z$.

Such a ring ${\cal O}$ has rationalization a quadratic extension of
$\mb Q$, but might not be integrally closed.  We could indeed choose
such subrings of ${\cal O}_F$, and these would give more general
theories with interesting content, but ${\cal O}_F$ is a legitimate
starting point.

If we did not specify $F$ or ${\cal O}$ as part of the data, they
would be invariants of connected components.

\begin{ques}
Why do we require an action on the abelian variety itself?  Why don't
we simply require an abelian variety with a specified $1$-dimensional
summand of its $p$-divisible group?
\end{ques}

The short answer is that it is based on our desire for the Shimura
stack $\Sh$ to actually have some content at height $n$.

Essentially, any height $n$ point of such a moduli will automatically
have an action of a ring ${\cal O}_F$ for some $F$, or possibly a
subring ${\cal O}$ as specified in the previous question.  More,
simply specifying that we have a $1$-dimensional summand of the
$p$-divisible group will give a tremendous abundance of path
components of the moduli as in section \ref{sec:level}.  Those path
components that cannot be rectified to have ${\cal O}$-actions for
some ${\cal O}$ will not have any height $n$ points.

\begin{ques}
Why don't we simply pick a connected component of the moduli, rather
than specifying a level structure and possibly ending up with several
connected components?
\end{ques}

One problem is that it is hard to know how much data is required to
reduce down to a particular connected component, and even when it is
known it is hard to state it.  This kind of data is often a question
about class groups.

Even then, the resulting moduli is no longer defined over $\mb Z_p$,
but instead usually defined over some algebraic extension.

\begin{ques}
Which choices of quadratic imaginary field and level structure data
determine interesting Shimura varieties?  How does the structure of
the spectrum $\TAF$ vary depending on these inputs?
What does the global geometry of these moduli look like (in
characteristic $0$ or characteristic $p$) at interesting chromatic
heights?  How does one go about computing these rings of integral,
or even rational, automorphic forms and higher cohomology?
\end{ques}

Some progress has been made at understanding chromatic level 2 and the
connection between $\TAF$ and $\TMF$.  A brief description of this is
to follow in section \ref{sec:absurf}.  The structure definitely
varies from input to input.  However, this is a place where more
computation is needed, and to create more computation one needs to use
more techniques for computing with these algebraic stacks that are not
simply presented by Hopf algebroids.

\section{Example: CM curves and abelian surfaces}
\label{sec:absurf}

We list here two basic examples of these moduli of abelian varieties
at chromatic levels $1$ and $2$.

At chromatic level $1$, the objects we are classifying are elliptic
curves with complex multiplication (the polarization data turns out to
be redundant).  Associated to a quadratic imaginary extension $F$ of
$\mb Q$, the moduli roughly takes the form
\[
\coprod_{{\rm Cl}(F)} [\ast // {\cal O}_F^\times].
\]
Here ${\rm Cl}(F)$ is the class group of $F$, and $[*//G]$ denotes a
point with automorphism group $G$.  This is, strictly speaking, only a
description of the geometric points of the stack.

At chromatic level $2$, the objects under study are abelian surfaces
with polarization and action of ${\cal O}_F$, together with a level
structure.  Ignoring the level structure, one can construct various
path components of the moduli as follows.  (This describes forthcoming
work \cite{absurf}.)

Given an elliptic curve $E$, we can form a new abelian surface $E
\otimes {\cal O}_F \cong E \times E$, with ${\cal O}_F$-action through
the second factor.  The Hermitian pairing on ${\cal O}_F$, together
with a ``canonical'' polarization on the elliptic curve $E$, gives
rise to a polarization of $E \otimes {\cal O}_F$ that
conjugate-commutes with the ${\cal O}_F$-action.  This construction is
natural in the elliptic curve, and produces a map of moduli
\[
{\cal M}_{ell} \to \Sh.
\]
The image turns out to be a path component of $\Sh$.  This is an
isomorphism onto the path component unless $F$ is formed by adjoining
a $4$th or $6$th root of unity.  In these cases it is a degree 2 or
degree 3 cover respectively, and we recover spectra with homotopy
\[
\mb Z_p[c_4,c_6^2,\Delta^{-1}] \subset \pi_* \TMF[\Delta^{-1}]
\]
for primes $p \equiv 1$ mod $4$, and
\[
\mb Z_p[c_4^3,c_6,\Delta^{-1}] \subset \pi_* \TMF[\Delta^{-1}].
\]
for primes $p \equiv 1$ mod $3$.

There are generalizations and modifications of this construction to
recover path components for other choices of level structure.  In
particular, by using alternate constructions we obtain objects which
are homotopy fixed points of the action of an Atkin-Lehner involution
on spectra $\TMF_0(N)[\Delta^{-1}]$.

Two such examples are as follows.  These rings of modular forms are
subrings of those described by Behrens \cite{marksphere} and
Mahowald-Rezk \cite{rezkmahowald} respectively.

If $p > 3$ is congruent to $1$ or $3$ mod $8$, there is a spectrum
associated to a moduli of abelian varieties with $\mb
Z[\sqrt{-2}]$-multiplication whose homotopy is a subring
\[
\mb Z_p[q_2, D^{\pm 1}]/ \subset
\TMF_0(2)[\Delta^{-1}]_*
\]
of the $p$-completed ring of modular forms of level $2$, where $|q_2|
= 4$ and $|D| = 8$.

If $p$ is congruent to $1$ mod $3$, there is a spectrum associated to
a moduli of abelian varieties with $\mb Z[({1 +
  \sqrt{-3}})/2]$-multiplication whose homotopy is a subring
\[
\mb Z_p[a_1^6, D^{\pm 1}]/ \subset
\TMF_0(3)[\Delta^{-1}]_*
\]
of the $p$-completed ring of modular forms of level $3$, where
$|a_1^6| = |D| = 12$.

\nocite{*}
\bibliography{tafoverview}

\end{document}